\documentclass[twoside,11pt]{amsart}

\usepackage{amsmath,latexsym,amssymb}%,times,mathptm}

\def\sqr#1#2{{\vcenter{\hrule height.#2pt
        \hbox{\vrule width.#2pt height#1pt \kern#1pt
                \vrule width.#2pt}
        \hrule height.#2pt}}}

\def\m{{\mathfrak m}}
\def\p{{\mathfrak p}}
\def\frak{\mathfrak}

\def\demo{\noindent{\bf Proof. }}
\def\QED{\hfill$\Box$}

\newcommand{\rar}{\rightarrow}
\newcommand{\lar}{\longrightarrow}

\newtheorem{Theorem}{Theorem}[section]
\newtheorem{Lemma}[Theorem]{Lemma}
\newtheorem{Corollary}[Theorem]{Corollary}
\newtheorem{Proposition}[Theorem]{Proposition}

\newtheorem{Question}[Theorem]{Question}
\theoremstyle{definition}

\theoremstyle{remark}

\newtheorem{Remark}[Theorem]{Remark}

\begin{document}

\baselineskip=17pt

\title[Normalization of Ideals and Brian\c{c}on--Skoda Numbers]
{\Large\bf Normalization of ideals and Brian\c{c}on--Skoda
numbers}

\author[C. Polini, B. Ulrich and W.V. Vasconcelos]
{Claudia Polini \and Bernd Ulrich \and Wolmer V. Vasconcelos}

\thanks{AMS 2000 {\em Mathematics Subject Classification}.
Primary 13A30; Secondary 13B22, 13H10, 13H15.}

\thanks{The three authors gratefully acknowledge partial support from the NSF}

\address{Department of Mathematics, University of Notre Dame,
Notre Dame, Indiana 46556} \email{cpolini@nd.edu}
\urladdr{www.nd.edu/{\textasciitilde}cpolini}

\address{Department of Mathematics, Purdue University,
West Lafayette, Indiana 47907} \email{ulrich@math.purdue.edu}
\urladdr{www.math.purdue.edu/{\textasciitilde}ulrich}

\address{Department of Mathematics, Rutgers University,
Piscataway, New Jersey 08854} \email{vasconce@math.rutgers.edu}
\urladdr{www.math.rutgers.edu/{\textasciitilde}vasconce}

\vspace{0.1in}

\begin{abstract} We establish bounds for the coefficient
$\overline{e}_1(I)$ of the Hilbert function of the integral closure
filtration of equimultiple ideals. These values are shown to help control
all algorithmic processes of normalization that make use of
extensions satisfying the condition $S_2$ of Serre.
\end{abstract}

\maketitle

\vspace{0.0in}

\section{Introduction}
Let $R$ be a Noetherian ring and let $I$ be an $R$--ideal. The
{\it integral closure} of $I$ is the ideal $\overline{I}$
consisting of all $z \in R$ which are solutions of equations of
the form
\[ z^n + a_1z^{n-1} + \, \ldots \, + a_n=0, \ \ \ a_i \in I^{i} .\]
The authors are not aware of any direct algorithm that builds
$\overline{I}$ from $I$, a situation that is aggravated by the
lack of numerical measures to distinguish between the two ideals.
A {\em non direct} construction of the integral closure of an
ideal passes through the integral closure $\overline{R[It]}$ in
$R[t]$ of the Rees algebra $R[It]$. In fact $\overline{I}$ is the
degree one component of $\overline{R[It]}$,
\[ I \leadsto \overline{R[It]} = R \oplus
\framebox{$\overline{I}t$} \oplus \overline{I^2}t^2 \oplus \,
\ldots \ \leadsto \overline{I}.\] This begs the issue since the
construction of $\overline{R[It]}$
%, for
%arbitrary ideals, may
%verge on the impossibility. It
 takes place in a much larger setting,
%(that of a presentation $R[It]= R[T_1, \ldots, T_n]/L$).
%By
while  in a {\em direct} construction $I\leadsto \overline{I}$
%we mean an
the steps of the algorithm %whose steps
would  take place entirely in $R$.
% These are lacking
%in the literature. We provide in full details a simple situation,
%that of monomial ideals of finite colength.

We refer to $\overline{R[It]}$ as the {\em normalization} of $I$.
Its construction is a standard step in the theory of
desingularization. It is very significant that there are numerical
measures to tell the two algebras $R[It]$ and $\overline{R[It]}$
apart for classes of ideals of interest. We are going to identify
one such measure, establish bounds for its value on
$\overline{R[It]}$ by data on $I$, and show how it bounds the
number of iterations of any algorithm that builds
$\overline{R[It]}$ by a succession of graded extensions
\[ R[It]  \rar A_1 \rar A_2 \rar \, \ldots \, \rar A_n =\overline{R[It]}\]
satisfying the condition $S_2$ of Serre. Recall that if one
chooses $A_1$ to be the $S_2$--ification of $R[It]$, then the
algorithm of \cite{clos} indeed produces such chains of $S_2$
algebras, provided $R$ is a $S_2$ domain of dimension $\geq 2$ and
is essentially of finite type over a field of characteristic zero.

\smallskip

We now describe the results of the paper. Let $(R, \mathfrak{m})$
be  an analytically unramified local Cohen--Macaulay ring of
dimension $d \geq 2$ and type $t$ with infinite residue field, $I$
an $\mathfrak{m}$--primary ideal, and $R[It] \subset A \subset B
\subset \overline{R[It]}$ inclusions of graded $R$--algebras. With
$e_1(A)$ denoting the first Hilbert coefficient of the ideal
filtration $\{A_n\}$ arising from $A$, we prove in Theorem 2.2
that $e_1(A) < e_1(B)$, provided $A$ satisfies $S_2$ and is
properly contained in $B$. The monotonicity of the function
$e_1(-)$ on these algebras yields the crucial role of
$\overline{e}_1(I)=e_1(\overline{R[It]})$ in a numerical criterion
of normality. When $R$ has a canonical module the $S_2$--ification
of $R[It]$ is given by ${\rm End}_{R[It]}(\omega_{R[It]})$, which
is relatively easy to compute. Now Corollary~\ref{normality} says
that $\overline{R[It]}$ can be obtain as the $S_2$--ification of
$R[It]$ if and only if $\overline{e}_1(I)=e_1(I)$; in particular
$I$ is normal if and only if $R[It]$ satisfies $S_2$ and
$\overline{e}_1(I)=e_1(I)$. In general $\overline{e}_1(I)$ can be
used to bound the lengths of chains of graded $S_2$ algebras lying
between $R[It]$ and $\overline{R[It]}$, see
Corollary~\ref{DivChains}.
%At this point one naturally asks

Thus one is led to search for effective upper bounds on
$\overline{e}_1(I)$. Notice that any such inequality also bounds
the first Hilbert coefficient $e_1(I)$, an issue that has been
addressed in \cite{HM, E, E2, RV} for instance. The bounds we are
looking for should estimate $\overline{e}_1(I)$ in terms of the
multiplicity $e_0(I)$ of the ideal $I$. The link between these two
Hilbert coefficients is provided by the {\it Brian\c{c}on--Skoda
number} $b(I)$ of $I$, which is the smallest integer $b$  such
that $\overline{I^{n+b}}\subset J^n$ for every $n$ and every
reduction $J$ of $I$. Indeed, in Theorem~\ref{maintheorem}(c) we
prove that
\[\overline{e}_1(I) \leq b(I) \ {\rm min} \, \{\frac{t}{t + 1} \
e_0(I), e_0(I) -\lambda(R/\overline{I})\},\] where $\lambda(-)$
denotes length. Furthermore, motivated by this result we estimate
the Brian\c{c}on--Skoda number of $I$ in
Proposition~\ref{boundonBS}. In a regular local ring the above
inequality reads \[ \overline{e}_1(I) \leq (d-1) \ {\rm min}\,
\{\frac{e_0(I)}{2}, e_0(I) -\lambda(R/\overline{I})\},\] since in
this case $b(I) \leq d-1$ by the classical Brian\c{c}on--Skoda
theorem. If in addition $d=2$ and $I$ is integrally closed we
obtain the equalities $\overline{e}_1(I)=e_1(I) = e_0(I) -
\lambda(R/I)$, which in turn imply the well--known facts that $I$
has reduction number at most one and $R[It]$ is Cohen--Macaulay
and normal.

In Section 3 we also establish bounds for $\overline{e}_1(I)$ that
avoid any reference to the Brian\c{c}on--Skoda number and instead
only involve the multiplicities of $I$ and of $I$ modulo an
element in the Jacobian of $R$. Our proofs are based on a general
Brian\c{c}on--Skoda type theorem due to Hochster and Huneke that
applies to non regular rings as well.
% in \cite[1.5.5 and 4.1.5]{HH}.
In Theorem~\ref{maintheorem}(a),(b) we show that if $R$ is an algebra
essentially of finite type over a perfect field $k$
and $\delta$ is a non zerodivisor in ${\rm Jac}_k(R)$, then
\[
\overline{e}_1(I) \leq \frac{t}{t + 1}\bigl[(d-1)e_0(I) + e_0(I
+\delta R/\delta R)\bigr]\] \,\, and
\[ \overline{e}_1(I) \leq (d-1)\bigl[e_0(I)
-\lambda(R/\overline{I})\bigr] + e_0(I +\delta R/\delta R).\] In
Section 4 we extend these results to arbitrary equimultiple
ideals.

\bigskip

\section{Normalization of Rees algebras}

The computation (and its control) of the integral closure of a
standard graded algebra over a field benefits greatly from Noether
normalizations and of the structures built upon them. If $A=R[It]$
is the Rees algebra of an ideal $I$ in a Noetherian ring $R$, it
does not allow for many such constructions. We would  still like
to develop some tracking of the complexity of the task required to
build $\overline{A}$ (assumed $A$--finite) through sequences of
graded extensions
\[ A \rar A_1 \rar A_2 \rar \cdots \rar A_n =\overline{A}\]
where $A_{i+1}$ is obtained from a specific procedure
%$\mathcal{P}$
applied to $A_i$. As in  \cite{teneadd}, if the
$A_i$ satisfy  the condition $S_2$ of Serre, we will call such
chains {\it divisorial}. At a minimum, we would want to bound the
length of divisorial chains. In this section we show how this can
be realized for Rees algebras of ideals.

\medskip

We now review some definitions and basic facts. For ideals $J
\subset I$ in a Noetherian ring one says that $J$ is a {\it
reduction} of $I$ if $\overline{I}=\overline{J}$. If in addition
$R$ is local with infinite residue field, we define {\it minimal}
reductions of $I$ to be reductions minimal with respect to
inclusion. The minimal number of generators of every minimal
reduction of $I$ is the {\it analytic spread} of $I$, which is
bounded below by the height of the ideal $I$ and above by the
dimension of the ring $R$. Thus if $I$ is an  $\m$--primary ideal
every minimal reduction of $I$ is generated by ${\rm dim} \, R$
elements. Finally, we say that $I$ is {\it equimultiple} if every
minimal reduction of $I$ is generated by ${\rm ht} \, I$ elements.

\medskip

Let $(R, \mathfrak{m})$ be  a Noetherian local ring of dimension
$d > 0$ and let $I$ be an $\mathfrak{m}$--primary ideal. Let $D=
\bigoplus_{n \geq 0} D_n t^n$ be a graded $R$--subalgebra of
$R[t]$ with $R[It] \subset D \subset R[t]$ and assume that $D$ is
a finite $R[It]$--module. For any such algebra we consider the
Hilbert--Samuel function $\lambda(R/D_n)$.
% where $\lambda(-)$ denotes length.
For $n\gg 0$ this function is given by the
Hilbert--Samuel polynomial
\[ e_0(D){{n+d-1}\choose{d}}
-e_1(D){{n+d-2}\choose{d-1}}+ \textrm{\rm lower terms}\,.\]
Notice that $e_i(R[It])$ coincide with the usual {\it Hilbert
coefficients} $e_i(I)$ of $I$. Furthermore $e_0(D) = e_0(I)$. By $
\overline{R[It]}$ we will always denote the integral closure of
$R[It]$ in $R[t]$. We write $\overline{e}_i(I)$ for the {\it
normalized Hilbert coefficients} $e_i(\overline{R[It]})$ of $I$ in
case $ \overline{R[It]}$ is a finite $R[It]$--module. The
% Recall that
%this assumption is always satisfied if $R$ is an analytically
%unramified local ring (\cite[]{R}). The
%first {\it normalized Hilbert
coefficient $\overline{e}_1(I)$ will be the main object of
interest in this paper.
%(Note to ourselves: the notation was chosen to distinguish it from
%$e_1(\overline{I})$.)

\medskip

%The above assumption that $ \overline{R[It]}$ be a finite
%$R[It]$--module holds for any ideal $I$ in an analytically
%unramified local ring $R$, see \cite[]{R}. Thus in this case there
%exits an integer $b$ such that
The condition that $ \overline{R[It]}$ be a finite $R[It]$--module
is satisfied for any ideal $I$ in an analytically unramified local
ring $R$, see \cite[1.5]{R}. Under this assumption there exits an
integer $b$ such that
%For arbitrary ideal $I$ the in a
%Noetherian ring $R$ we define the {\it Brian\c{c}on--Skoda number}
%$b(I)$ of $I$ to be the smallest integer $b\geq 0$ such that
\[
\overline{I^{n+b}} \subset J^n \quad {\rm for \ every \ } n \ {\rm
and \ every \ reduction \ } J \   {\rm of \ } I,
\]
where we use the convention that $I^m=R$ for $m\leq 0$.  The
smallest such $b\geq 0$ is called the {\it Brian\c{c}on--Skoda
number} $b(I)$ of $I$.
%If in addition $R$ is excellent then by
%\cite[]{Hu}, the Brian\c{c}on--Skoda number of any ideal $I$ is
%uniformly bounded by an integer $c$ only depending on $R$. One can
%then define $c(R)$ to be the smallest such $c$.
In a regular local ring $R$ of positive dimension one has
$b(I)\leq\dim R -1$ according to the classical Brian\c{c}on--Skoda
theorem (\cite[Theorem 1]{LipmanSathaye}).

\medskip

In order to relate  $e_0(I)$, $\overline{e}_1(I)$ and $b(I)$ the
next lemma is needed. We use the notation $\deg(-)$ for the
multiplicity of a finite module over a Noetherian local ring or of
a finite graded module over a Noetherian standard graded algebra
over an Artinian local ring.

\begin{Lemma}\label{lemma1}
Let $(R, \mathfrak{m})$ be a Noetherian local ring of dimension
$d>0$ and let $I$ be an $\mathfrak{m}$--primary ideal. Let $A$ and
$B$ be graded $R$--subalgebras of $R[t]$ with
\[ R[It] \subset A \subset B \subset R[t]\]
and assume that $B$ is a finite $R[It]$--module. Write $C$ for the
graded $R[It]$--module $B/A$.

\begin{itemize}
\item[(a)] ${\rm dim} \, C \leq d$, and equality holds if $R$ is
Cohen--Macaulay, $A$ satisfies the condition $S_2$ of Serre and
$A\not=B$.

\item[(b)] If ${\rm dim} \, C < d$ then $e_1(B)=e_1(A)$.

\item[(c)] If ${\rm dim} \, C = d$ then $e_1(B)-e_1(A)= \deg(C) >
0$.
\end{itemize}
\end{Lemma}
\demo Part (a) is obvious. To prove (b) and (c) notice that $C$ is
a finite graded module over a Noetherian standard graded algebra
over an Artinian ring. Hence it has a Hilbert polynomial whose
degree is ${\rm dim} \, C -1$. On the other hand
%${\rm dim} \, A={\rm dim} \, B=d+1$. Hence
the exact sequences
\[ 0 \lar C_n \lar R/A_n \lar R/B_n \lar 0\]
show that $(e_1(B)-e_1(A))/(d-1)!$ is the coefficient of the term
of degree $d-1$ in the Hilbert polynomial of $C$.\QED

\medskip

\begin{Theorem}\label{e1versusb}
Let $(R, \mathfrak{m})$ be  an analytically unramified local
Cohen--Macaulay ring of positive dimension with infinite residue
field and let $I$ be an $\mathfrak{m}$--primary ideal. Let $A$ and
$B$ be distinct graded $R$--subalgebras of $R[t]$ with
\[ R[It] \subset A \subsetneq B \subset \overline{R[It]}\]
and assume that $A$ satisfies the condition $S_2$ of Serre.
Then
\[ 0 \leq e_1(I) \leq e_1(A) <  e_1(B) \leq \overline{e}_1(I) \leq b(I)e_0(I).\]
\end{Theorem}
\demo  Let $J$ be a minimal reduction of $I$, and notice that
$e_1(J)=0$ since $R$ is Cohen--Macaulay. Now Lemma~\ref{lemma1}
implies the asserted inequalities except for the last one. To show
that $\overline{e}_1(I) \leq b(I)e_0(I)$ write $d= {\rm dim} \, R$
and $b=b(I)$. We consider the $R[Jt]$--module $C=\overline{R[
It]}/R[Jt]$.  Since $R[Jt]$ satisfies the condition $S_2$ of
Serre, Lemma~\ref{lemma1} shows that \[\overline{e}_1(I)=
\overline{e}_1(I)-e_1(J) = \deg(C).\] In turn, by the definition
of $b$, $C$ is a submodule of the graded $R[Jt]$--module $D$ whose
$n$--th components are $J^{n-b}/J^{n}$.
%As ${\rm dim}\, D=d={\rm
%dim}\, C$ one has $\deg(C) \leq \deg(D)$. To determine the latter
The inclusions $J^{n-b} \supset J^{n-b+1} \supset \ldots \supset
J^{n}$ induce a filtration of $D$. From their Hilbert functions
one sees that the factors in this filtration all have dimension
$d$ and multiplicity $e_0(J)=e_0(I)$. Hence ${\rm dim}\, D=d={\rm
dim}\, C$ and $\deg(D)=b \ e_0(I)$. Now the containment $C \subset
D$ gives
\[ \deg(C) \leq \deg(D)= b \ e_0(I),\]
%Indeed, for each positive integer $k$,
%\[ \lambda(J^{n+k-1}/J^{n+k})= e_0(J){{n+k-1+d-1}\choose{d-1}}=
%\frac{e_0(J)}{(d-1)!} n^{d-1} + \textrm{\rm lower terms}.\]
and we obtain $\overline{e}_1(I) = \deg(C) \leq b(I)e_0(I)$, as
asserted.
 \QED

\medskip

\begin{Remark}\label{REMARK1}
{\rm Applying Theorem~\ref{e1versusb} to any minimal reduction $J$
of $I$ one obtains the stronger estimate $\overline{e}_1(I) =
\overline{e}_1(J) \leq  b(J)e_0(J) = b(J)e_0(I)$. }
\end{Remark}

\smallskip

A measure for the complexity of $\overline{R[It]}$ is the number
of steps needed to construct it. We address this issue in the next
corollary, which is a direct consequence of
Theorem~\ref{e1versusb}.

\begin{Corollary}\label{DivChains}Let $(R, \mathfrak{m})$ be
an analytically unramified local Cohen--Macaulay ring of positive
dimension and let $I$ be an $\mathfrak{m}$--primary ideal. Then
$\overline{e}_1(I)$ bounds the length of any chain of graded
$R$--subalgebras satisfying the condition $S_2$ of Serre lying
strictly between $R[It]$ and $\overline{R[It]}$.
\end{Corollary}

\smallskip

The following corollary provides a numerical criterion for when
the integral closure coincides with the $S_2$--ification ${\rm
End}_{R[It]}(\omega_{R[It]})$ of $R[It]$. Let $(R, \mathfrak{m})$
be a local Cohen--Macaulay ring of dimension $\geq 2$ with a
canonical module and infinite residue field, let $I$ be an
$\mathfrak{m}$--primary ideal, and $J$ a minimal reduction of $I$.
Writing $(-)^{\vee}={\rm Hom}_{R[Jt]}(-,R[Jt])$ we consider these
embeddings of graded algebras,

\[ R[It] \subset {\rm End}_{R[It]}(\omega_{R[It]}) \subset
{\rm End}_{R[It]}(R[It]^{\vee})=R[It]^{\vee\vee} \subset
\overline{R[It]}.\]
\smallskip
\begin{Corollary}\label{normality} Let $(R, \mathfrak{m})$ be an analytically unramified
local Cohen--Macaulay ring of dimension $\geq 2$ with a canonical
module and infinite residue field, and let $I$ be an
$\mathfrak{m}$--primary ideal. Then $\overline{R[It]}= {\rm
End}_{R[It]}(\omega_{R[It]})$ if and only if
$\overline{e}_1(I)=e_1(I)$$;$ in this case
$\overline{R[It]}=R[It]^{\vee\vee}$. In particular, the ideal $I$
is normal if and only if $R[It]$ satisfies the condition $S_2$ of
Serre  and $\overline{e}_1(I)=e_1(I)$.
\end{Corollary}
\demo Notice that by Lemma~\ref{lemma1}(b), ${\rm
End}_{R[It]}(\omega_{R[It]})$ has first Hilbert coefficient
$e_1(I)$. Now apply Theorem 2.2. \QED

\bigskip

\bigskip

\section{Bounds on $\overline{e}_1(I)$ and the Brian\c{c}on--Skoda number}

We discuss the role of Brian\c{c}on--Skoda type theorems (see
\cite{AHu}, \cite{LipmanSathaye}) in determining some
relationships between the coefficients $e_0(I)$ and
$\overline{e}_1(I)$. We will use a Brian\c{c}on--Skoda theorem
that works in non--regular rings. We are going to provide a short
proof along the lines of \cite{LipmanSathaye} for the special case
we need: $\mathfrak{m}$--primary ideals in a local Cohen--Macaulay
ring. The general case is treated by Hochster and Huneke in
\cite[1.5.5 and 4.1.5]{HH}. Let $k$ be a perfect field, let $R$ be
a reduced and equidimensional  $k$--algebra essentially of finite
type, and assume that $R$ is affine with $d={\rm dim} \, R$ or
$(R, \mathfrak{m})$ is local with $d={\rm dim}\, R + {\rm trdeg}_k
R/{\m}$. Recall that the {\it Jacobian ideal}   ${\rm Jac}_k(R)$
of $R$ is defined as the $d$--th Fitting ideal of the module of
differentials $\Omega_k(R)$ -- it can be computed explicitly from
a presentation of the algebra. By varying Noether normalizations
one deduces from \cite[Theorem 2]{LipmanSathaye} that the Jacobian
ideal ${\rm Jac}_k(R)$ is contained in the conductor $R \colon
\overline{R}$ of $R$ (see also \cite{N}, \cite[3.1]{AB} and
\cite[2.1]{H}); here $\overline{R}$ denotes the integral closure
of $R$ in its total ring of fractions.

\smallskip

\begin{Theorem}\label{GBS}
Let $k$ be a perfect field, let $R$ be a reduced local
Cohen--Macaulay $k$--algebra essentially of finite type, and let
$I$ be an equimultiple ideal of height $g > 0$.  Then for every
integer $n$,
\[
{\rm Jac}_k(R) \, \overline{I^{n+g -1}} \subset I^n.
\]
\end{Theorem}
\demo We may assume that $k$ is infinite. Then, passing to a
minimal reduction, we may suppose that $I$ is generated by a
regular sequence of length $g$. Let $S$ be a finitely generated
$k$--subalgebra of $R$ so that $R=S_{\p}$ for some $\p \in {\rm
Spec}(S)$, and write $S=k[x_1, \ldots, x_e] = k[X_1, \ldots,
X_e]/{\frak a}$ with ${\frak a}=(h_1, \ldots, h_t)$ an ideal of
height $c$. Notice that $S$ is reduced and equidimensional. Let
$K= (f_1, \ldots, f_g)$ be an $S$--ideal with $K_{\p}=I$, and
consider the extended Rees ring $B=S[Kt, t^{-1}]$. Now $B$ is a
reduced and equidimensional affine $k$--algebra of dimension
$e-c+1$.

Let $\varphi \colon k[X_1, \ldots, X_e, T_1, \ldots, T_g, U]
\twoheadrightarrow B$ be the $k$--epimorphism mapping $X_i$ to
$x_i$, $T_i$ to $f_it$ and $U$ to $t^{-1}$. Its kernel has height
$c+g$ and contains the ideal ${\frak b}$ generated by $\{ h_i,
T_jU-f_j | 1 \leq i \leq t, 1 \leq j \leq g \}$. Consider the
Jacobian matrix of these generators,
\[
\Theta = \left(
\begin{array}{c|cccc}
\displaystyle\frac{\partial h_i}{\partial X_j} & &   \!\!\!\! 0 & & \\
\hline
& U & & & T_1 \\
* & & \ddots & & \vdots \\
& & & U & T_g
\end{array}
\right).
\]
Notice that $ I_{c+g}(\Theta)\supset I_c(\left(\frac{\partial
h_i}{\partial X_j}\right))U^{g-1}(T_1, \ldots, T_g)$. Applying
$\varphi$ we obtain $ {\rm Jac}_k(B) \supset I_{c+g}(\Theta) B
\supset {\rm Jac}_k(S) K t^{-g+2}$. Thus ${\rm Jac}_k(S) K
t^{-g+2}$ is contained in the conductor of $B$. Localizing at $\p$
we see that ${\rm Jac}_k(R) I t^{-g+2}$ is in the conductor of the
extended Rees ring $R[It, t^{-1}]$. Hence for every $n$, ${\rm
Jac}_k(R) \, I \, \overline{I^{n+g-1}} \subset I^{n+1}$, which
yields
\[
{\rm Jac}_k(R) \, \overline{I^{n+g-1}} \subset I^{n+1} \colon I =
I^n,
\]
as $({\rm gr}_I(R))_{+}$ has positive grade. \QED

\bigskip

We now use Theorem~\ref{GBS} to sharpen the bound on
$\overline{e}_1$ given in Theorem~\ref{e1versusb}.

\begin{Theorem}\label{maintheorem} Let $(R, \m)$ be a
reduced local Cohen--Macaulay ring of dimension $d
> 0$ and let $I$ be an $\mathfrak{m}$--primary ideal.
\begin{itemize}
\item[$($a$)$] If in addition $R$ is an algebra essentially of
finite type over a perfect field $k$ with type $t$, and $\delta
\in {\rm Jac}_k(R)$ is a non zerodivisor, then \[\overline{e}_1(I)
\leq \frac{t}{t + 1} \bigl[ (d-1)e_0(I) + e_0(I +\delta R/\delta
R) \bigr].\]
 \item[$($b$)$] If the assumptions of $($a$)$ hold,
then
\[
\overline{e}_1(I) \leq (d-1)\bigl[e_0(I)
-\lambda(R/\overline{I})\bigr] + e_0(I +\delta R/\delta R).
\]
\item[$($c$)$] If $R$ is analytically unramified and $R/\m$ is
infinite, then
\[
\overline{e}_1(I) \leq b(I) \ {\rm min} \, \{\frac{t}{t + 1} \ e_0(I),
e_0(I) -\lambda(R/\overline{I})\}.
\]
\end{itemize}
\end{Theorem}
\demo We may assume that $R/\m$ is infinite. Then, passing to a
minimal reduction we may suppose that $I$ is generated by a
regular sequence $f_1, \ldots, f_d$. Notice this can only decrease
$b(I)$. Let $S$ be a local ring obtained from $R$ by a purely
transcendental residue field extension and by factoring out $d-1$
generic elements $a_1, \ldots, a_{d-1}$ of $I$. To be more
precise, $S=R(\{X_{ij}\})/(a_1, \dots, a_{d-1})$ with $\{X_{ij}\}$
a set of $(d-1)d$ indeterminates and
$a_i=\sum_{j=1}^{d}X_{ij}f_j$. Notice that $S$ is also a
birational extension of a localization of a polynomial ring over
$R$, and hence is analytically unramified according to
\cite[36.8]{Nag} and \cite[1.6]{R}. Furthermore $S$ is a
one--dimensional local ring and the $S$--ideal $IS$ is generated
by a single non zerodivisor, say $IS=fS$. From \cite[Theorem
1]{Itoh} one has
\begin{equation}\label{eq2}
\overline{I}S = \overline{IS},
\end{equation}
\begin{equation}\label{eq3}
\overline{I^{n}}S = \overline{(IS)^{n}} \quad {\rm for \ } n \gg
0.
\end{equation}
The last fact combined with the genericity of $a_1, \ldots,
a_{d-1}$ yields $ \overline{e}_1(I) =\overline{e}_1(IS)$. Moreover
$\overline{e}_1(IS)= \lambda(\overline{S}/S)$ as $S$ is a
one--dimensional analytically unramified local ring. Thus
\begin{equation}\label{eq4}
\overline{e}_1(I) =\lambda(\overline{S}/S).
\end{equation}

In the setting of $($a$)$ and $($b$)$ the element $\delta$ is a
non zerodivisor on $S$. Furthermore Theorem~\ref{GBS} shows that
\begin{equation}\label{eq1}
\delta I^{d-1} \subset \delta \ \overline{I^{d-1}}\subset I^n
\colon  \overline{I^{n}}  \qquad {\rm for \ every \ } n.
\end{equation}
For $n \gg 0$, by $($\ref{eq3}$)$, $($\ref{eq1}$)$ and since $f^n
S$ is contained in the conductor $S \colon \overline{S}$, we
obtain
\[
\delta f^{d-1} S \subset f^n S \colon \overline{f^n S} = f^n S
\colon f^n \overline{ S}=S \colon \overline{S}.
\]
Hence
\begin{equation}\label{eq5}
\delta f^{d-1} \overline{S} \subset S \colon \overline{S}.
\end{equation}

We prove $($a$)$ by computing lengths along the inclusions
\begin{equation}\label{eq6}
\delta f^{d-1} S \subset \delta f^{d-1} \overline{S} \subset S
\colon \overline{S} \subset S.
\end{equation}
Also recall that
\begin{equation}\label{eq7}
\lambda(\overline{S}/S) \leq t \ \lambda(S/S \colon \overline{S})
\end{equation}
by \cite[the proof of 3.6]{HK} (see also \cite[Theorem 1]{BrH} and \cite[2.1]{D}). We obtain
\begin{eqnarray*}
\frac{t + 1}{t} \  \overline{e}_1(I) & = & \lambda(\overline{S}/S)
+ \frac{1}{t} \ \lambda(\overline{S}/S) \qquad \ \ \ \ \qquad
\qquad \qquad  {\rm by}
\ (\ref{eq4})\\
& \leq &  \lambda(\delta f^{d-1}\overline{S}/\delta f^{d-1} S) +
\lambda(S/S \colon \overline{S}) \qquad \qquad {\rm by}
\ (\ref{eq7}) \\
& \leq  & \lambda(S/\delta f^{d-1} S)  \qquad \qquad \ \qquad
\qquad \qquad \qquad {\rm by} \
(\ref{eq6}) \\
& =  & (d-1) \ \lambda(S/f S) + \lambda(S/\delta S) \\
& =  & (d-1)e_0(I) + e_0(I +\delta R/\delta R)  \qquad \qquad {\rm
by \ the \ genericity \ of \ } a_1, \ldots, a_{d-1}.
\end{eqnarray*}

Next we prove part $($b$)$. The inclusion (\ref{eq5}) yields the
filtration
\[
%\begin{equation}
%\label{eq8}
S= \delta f^{d-1}\overline{S} + S \subset f^{d-1}\overline{S} + S
\subset \ldots \subset f^2\overline{S} + S \subset f \overline{S}
+ S \subset \overline{S}, \]
%\end{equation}
which shows
\begin{equation}\label{eq9}
 \overline{e}_1(I) =\lambda(\overline{S}/S)=  \sum_{i=1}^{d-1}
 \lambda(f^{i-1}\overline{S} + S/ f^i\overline{S} + S) +
\lambda(f^{d-1}\overline{S} + S/ \delta f^{d-1} \overline{S} + S).
\end{equation}
Multiplication by $f$ induces epimorphisms of $S$--modules
\begin{equation}\label{eq10}
\frac{f^{i-1}\overline{S} + S}{ f^i\overline{S} + S}
\twoheadrightarrow \frac{f^{i}\overline{S} + S}{
f^{i+1}\overline{S} + S}\ .
\end{equation}
Now (\ref{eq9}) and (\ref{eq10}) show
\begin{equation}\label{eq11}
 \overline{e}_1(I) \leq (d-1) \ \lambda(\overline{S}/f\overline{S} +S) +
\lambda(f^{d-1}\overline{S} + S/ \delta f^{d-1} \overline{S} + S
).
\end{equation}
Next we claim that
\begin{equation}\label{eq12}
\lambda(\overline{S}/f\overline{S} +S) = \lambda(\overline{I}/I).
\end{equation}
Indeed,
\begin{eqnarray*}
\lambda(\overline{S}/f\overline{S} +S)  & = &
\lambda(\overline{S}/f\overline{S}) -
\lambda(f\overline{S} +S/f\overline{S}) \\
& =&  \lambda(S/f S) -  \lambda(S/S \cap f \overline{S}) \\
& =  & \lambda(S \cap f \overline{S}/f S)  \\
& =  & \lambda(\overline{f S}/f S)\\
& =  &  \lambda(\overline{I}/I) \qquad \qquad {\rm by} \
(\ref{eq2}) .
\end{eqnarray*}
On the other hand,
\begin{eqnarray*}
\lambda(f^{d-1}\overline{S} + S/ \delta f^{d-1} \overline{S} + S )
&\leq&
\lambda(f^{d-1}\overline{S} / \delta f^{d-1} \overline{S}) \\
& =  &   \lambda(\overline{S} / \delta \overline{S}) \\
& =  &   \lambda (S / \delta S)\\
& =  &   e_0(I +\delta R/\delta R) \qquad \qquad {\rm by \ the \
genericity \ of \ } a_1, \ldots, a_{d-1}.
\end{eqnarray*}
Therefore
\begin{equation}\label{eq13}
\lambda(f^{d-1}\overline{S} + S/ \delta f^{d-1} \overline{S} + S )
\leq
 e_0(I +\delta R/\delta R).
\end{equation}
 Combining (\ref{eq11}),  (\ref{eq12}) and (\ref{eq13}) we deduce
\begin{eqnarray*}
 \overline{e}_1(I) &\leq&  (d-1) \ \lambda(\overline{I}/I) + e_0(I +\delta R/\delta R)
 \\
& =  &   (d-1)\bigl[e_0(I) -\lambda(R/\overline{I})\bigr] + e_0(I
+\delta R/\delta R).
\end{eqnarray*}

Finally we prove part $($c$)$.  Write $b=b(I)$. We first claim
that
\begin{equation}\label{eq14}
 f^{b} \overline{S} \subset S \colon \overline{S}.
\end{equation}
Indeed, for $n \gg 0$
\begin{eqnarray*}
f^n S & \supset &
\overline{I^{n + b}}S   \\
& = & \overline{(IS)^{n + b}}  \qquad \qquad {\rm by} \
(\ref{eq3}) \\
& =  &  \overline{f^{n +b} S}  \\
& =  & f^{n +b} \overline{S}   \qquad \qquad {\rm since \ } n \gg
0.
\end{eqnarray*}
Therefore $f^{b} \overline{S} \subset S $, proving (\ref{eq14}).
Now (\ref{eq14}) yields the filtrations
\begin{equation}\label{eq15}
 f^{b} S \subset f^b  \overline{S} \subset S
\colon \overline{S} \subset S,
\end{equation}
\begin{equation}\label{eq16}
S= f^b\overline{S} + S \subset \ldots \subset f^2\overline{S} + S
\subset f \overline{S} + S \subset \overline{S}.
\end{equation}
Filtration (\ref{eq15}) implies
\begin{eqnarray*}
\frac{t + 1}{t} \  \overline{e}_1(I) & = & \lambda(\overline{S}/S)
+ \frac{1}{t} \ \lambda(\overline{S}/S)  \ \ \qquad \qquad {\rm
by}
\ (\ref{eq4})\\
& \leq &  \lambda(f^{b}\overline{S}/f^{b} S) + \lambda(S/S \colon
\overline{S})   \ \qquad {\rm by}
\ (\ref{eq7}) \\
& \leq  & \lambda(S/f^{b} S)  \ \ \  \qquad \qquad \qquad \qquad
{\rm by} \
(\ref{eq15}) \\
& =  & b \ \lambda(S/f S) \\
& =  & b \ e_0(I) \qquad \qquad {\rm by \ the \ genericity \ of \
} a_1, \ldots, a_{d-1}.
\end{eqnarray*}
On the other hand filtration (\ref{eq16}) yields
\begin{eqnarray*}
\overline{e}_1(I)  & = & \lambda(\overline{S}/S) \\
& =&   \sum_{i=1}^{b}
 \lambda(f^{i-1}\overline{S} + S/ f^i\overline{S} + S) \\
& \leq & b \ \lambda(\overline{S}/f\overline{S} +S)   \qquad
\qquad {\rm by} \
(\ref{eq10})\\
& =  & b \ \lambda(\overline{I}/I) \qquad \qquad \qquad \ \  {\rm
by} \
(\ref{eq12}) \\
& =  &  b  \ (e_0(I)
 -\lambda(R/\overline{I})).
\end{eqnarray*}
\QED

\bigskip

\begin{Remark}
{\rm The multiplicity $e_0(I+\delta R/\delta R)$ occurring in
Theorem~\ref{maintheorem} can be bounded by
\[
e_0(I+\delta R/\delta R)\leq (d-1)! \  e_0(I) \, \deg(R/\delta R),
\]
where $\deg(R/\delta R)$ is the multiplicity of the local ring
$R/\delta R$. Indeed, \cite[Theorem 3]{Lech} gives
\[
e_0(I+\delta R/\delta R)\leq (d-1)! \  \lambda(R/I + \delta R) \,
\deg(R/\delta R) .\] }
\end{Remark}

\medskip

\begin{Corollary}\label{RLR} Let $(R, \mathfrak{m})$
be a regular local ring of dimension $d>0$ and let $I$ be an
$\mathfrak{m}$--primary ideal. Then
\[
e_1(I) \leq \overline{e}_1(I) \leq (d-1) \ {\rm
min}\{\frac{e_0(I)}{2}, e_0(I) -\lambda(R/\overline{I})\}.
\]
\end{Corollary}
\demo We may assume that $R/\m$ is infinite. The classical
Brian\c{c}on--Skoda theorem gives that $b(I) \leq d-1$, see
\cite[Theorem 1]{LipmanSathaye}. The assertions now follow from
Theorems~\ref{e1versusb} and \ref{maintheorem}$($c$)$. \QED

\bigskip

We are now going to use Corollary~\ref{RLR} to bound the length of
divisorial chains for  classes of Rees algebras.

\begin{Corollary} Let $(R, \mathfrak{m})$ be a regular local ring of
dimension $d > 0$ and let $I$ be an $\mathfrak{m}$--primary ideal.
Then $(d-1) \ {\rm min}\{\frac{e_0(I)}{2}, e_0(I)
-\lambda(R/\overline{I})\}$ bounds the length of any chain of
graded $R$--subalgebras satisfying the condition $S_2$ of Serre
lying strictly between $R[It]$ and $\overline{R[It]}$.
\end{Corollary}
\demo  The assertion follows from Corollaries~\ref{DivChains} and
\ref{RLR}. \QED

\medskip

\begin{Remark}
{\rm \begin{itemize}

\item[$($a$)$] Let $(R, \mathfrak{m})$ be a regular local ring of
dimension $2$ and let $I$ be an $\mathfrak{m}$--primary integrally
closed ideal. Then $ e_1(I)  = \overline{e}_1(I) = e_0(I)
-\lambda(R/I)\leq \frac{e_0(I)}{2}$. This follows, for instance,
from Corollary~\ref{RLR} combined with the inequality $e_1(I) \geq
e_0(I) - \lambda(R/I)$, see \cite[Theorem 1]{No}. Furthermore by
\cite[2.1]{Hu2} or \cite[3.3]{O}, the equality $e_1(I) = e_0(I) -
\lambda(R/I)$ implies that $I$ has reduction number at most one if
$R/\mathfrak{m}$ is infinite, a fact proved in \cite[5.4]{LiTe}.
Thus $R[It]$ is Cohen--Macaulay according to \cite[3.1]{valval}
and \cite[3.10]{GS}. Now Corollary~\ref{normality} yields the
well--known result that $I$ is normal, see \cite[Theorem $2^{'}$,
p. 385]{Z}.

\item[$($b$)$] Let $k$ be an infinite field, write $R=k[X_1,
\ldots, X_d]_{(X_1, \ldots, X_d)}$, let $\mathfrak{m}$ denote the
maximal ideal of $R$, and let $I$ be an $\mathfrak{m}$--primary
$R$--ideal generated by homogeneous polynomials in $k[X_1, \ldots,
X_d]$ of degree $s$. Then $ \overline{e}_1(I) =
e_1(\mathfrak{m}^s)= \frac{d-1}{2} e_0(I) (1- \frac{1}{s})\approx
\frac{d-1}{2} e_0(I)$. This shows that the estimate of
Corollary~\ref{RLR} is essentially sharp.
\end{itemize}}
\end{Remark}
\medskip

\begin{Proposition}\label{boundonBS} Let $k$ be a perfect
field, let $(R, \mathfrak{m})$  be a
reduced local Cohen--Macaulay $k$--algebra essentially of finite
type of dimension $d > 0$, and let $\delta \in {\rm Jac}_k(R)$ be
a non zerodivisor. Then for any $\mathfrak{m}$--primary ideal $I$,
\[
b(I) \leq  d-1 + e_0(I +\delta R / \delta R).
\]
\end{Proposition}
\demo We may assume that $R/\m$ is infinite. Then, replacing $I$
by a minimal reduction with the same Brian\c{c}on--Skoda number we
may suppose that $I$ is generated by a regular sequence of length
$d$. As in the proof of Theorem~\ref{maintheorem} let $S$ be a
local ring obtained from $R$ by a purely transcendental residue
field extension and by factoring out $d-1$ generic elements $a_1,
\ldots, a_{d-1}$ of $I$. Write $IS=fS$ and let $b$ be the smallest
non negative integer with $f^b \overline{S} \subset S \colon
\overline{S}$.
%\[
%b={\rm min} \, \{ \, n \geq 0 \ | \ f^n \overline{S} \subset S
%\colon \overline{S} \ \}.
%\]

We first claim that
\begin{equation}\label{eq17}
b(I) \leq b.
\end{equation}
Indeed, for any integer $n \geq 0$ we have
\[
\overline{I^{n + b }} S \subset \overline{I^{n + b} S} \subset
f^{n + b }\overline{S} \subset f^n S,
\]
hence $\overline{I^{n+b}}S \subset I^nS$.
%(I^n, a_1, \ldots, a_{d-1}).$
As ${\rm gr}_I(R)$ is a polynomial ring in $d$ variables over
$R/I$, the generic choice of $a_1, \ldots, a_{d-1}$ gives that
%the images of these variables in ${\rm gr}_{IS}(S)$ are still
%algebraically independent over $R/I$. Thus
${\rm gr}_I(R)$ embeds into ${\rm gr}_{IS}(S)$. Therefore
$\overline{I^{n+b}} \subset I^{n}$, proving (\ref{eq17}).

By (\ref{eq17}) it suffices to show that $b \leq d-1 + e_0(I +
\delta R/\delta R)$. To this end we may assume $b \geq d-1$. The
definition of $b$ yields the filtration

\begin{equation}\label{eq18}
S= f^b\overline{S} + S \subset \ldots \subset f^{d-1}\overline{S}
+ S.
\end{equation}
On the other hand (\ref{eq5}) implies
\begin{equation}\label{eq19}
S= \delta f^{d-1}\overline{S} + S \subset f^{d-1}\overline{S} + S.
\end{equation}
If $f^{i}\overline{S} + S = f^{i-1}\overline{S} + S$ for some $b
\geq i \geq d$, then multiplication by $f^{b-i}$ yields
$S=f^b\overline{S} + S=f^{b-1}\overline{S} + S$, contradicting the
minimality of $b$. Thus (\ref{eq18}) gives
\[
\lambda(f^{d-1}\overline{S} + S/S) \geq b-d+1.
\]
On the other hand from (\ref{eq13}) and (\ref{eq19}) we deduce
\[
\lambda(f^{d-1}\overline{S} + S/S) \leq e_0(I + \delta R/\delta
R).
\]
Thus $b -d +1 \leq e_0(I + \delta R/\delta R).$ \QED

\bigskip

\begin{Remark}
{\rm In the setting of the proof of Proposition~\ref{boundonBS},
$b(I)=b(IS)=b$. Indeed,  (\ref{eq3}) implies that for $n \gg 0$,
\[
f^{n + b(I)} \overline{S} = \overline{I^{n + b(I)} S} =
\overline{I^{n + b(I)}} S \subset I^n S = f^n S,
\]
showing that $f^{b(I)} \overline{S} \subset S.$ Hence $b \leq
b(I)$ and then $b =b(I)$ by (\ref{eq17}). Clearly $b(IS) \leq b$.
For $n \gg 0$, $f^{n + b(IS)}\overline{S} = \overline{f^{n +
b(IS)}S} \subset f^n S$ and therefore $f^{b(IS)} \overline{S}
\subset S$, showing that $b \leq b(IS)$.}
\end{Remark}

\bigskip

\section{Equimultiple ideals}

In this section we extend Theorems~\ref{e1versusb} and
\ref{maintheorem} to arbitrary equimultiple ideals.
%We place
%ourselves in the setting of these results with $I$ an equimultiple
%ideal of height $g$.
%Without harm we may assume that $I$ is a regular sequence.
The technical change involves Hilbert functions. Let $R$ be a
Noetherian local ring with infinite residue field and $I$ an ideal
of height $g>0$. Let $D= \bigoplus_{n \geq 0} D_n t^n$ be a graded
$R$--subalgebra of $R[t]$ with $R[It] \subset D \subset R[t]$ and
assume that $D$ is a finite $R[It]$--module. Instead of the length
function as in Theorem~\ref{e1versusb}, we consider the
multiplicity $\deg(R/D_n)$ of the $R$--module $R/D_n$. According
to the associativity formula for multiplicities one has
%, this function is
%easy to express in terms of localizations of $D_n$ at the minimal
%primes of $I$,
\[ \deg(R/D_n)= \sum_{\mathfrak p} \lambda(R_{\mathfrak p}/({D_n})_{\mathfrak
p})\deg(R/{\mathfrak p}),\] where $n \gg 0$ and the sum is taken
over the minimal primes $\mathfrak p$ of $I$ with $\dim
R/{\mathfrak p}=\dim R/I$. It follows that this function behaves
as a polynomial of degree $g$, which we still call the
Hilbert--Samuel polynomial of $D$,
\[  E_0(D){{n+g-1}\choose{g}}
-E_1(D){{n+g-2}\choose{g-1}}+ \textrm{\rm lower terms}\, .\] The
coefficients  $E_i(D)$ can be expressed in terms of the local
Hilbert coefficients $e_i(D_{\mathfrak{p}})$,
\[ E_i(D) = \sum_{\mathfrak p} e_i(D_{\mathfrak p})\deg(R/{\mathfrak p}).\]
If $g=0$ we set $E_0(D)=\sum_{\mathfrak p} \lambda(R_{\mathfrak
p}) \deg(R/{\mathfrak p})$. We will use the notation $E_i(I)$ when
$D=R[It]$ and $\overline{E}_i(I)$ when $D= \overline{R[It]}$, the
integral closure  of $R[It]$ in $R[t]$.

\smallskip
%\begin{Remark}{\rm
One glaring difficulty with the above formula lies in the numbers
$e_i(D_{\mathfrak p})$ or even $e_i(I_{\mathfrak p})$, which are
hard to get hold of. At least for equimultiple ideals on the other
hand, the qualitative behavior of the $E_i(I)$ is that of the
usual Hilbert coefficients and $E_0(I)$ can be expressed as a
multiplicity:

%%%precisely when $I$ is generically a complete intersection.
%}\end{Remark}

%There is a special case when some of these coefficients may be
%%infinite residue field and suppose $I$ is {\it equimultiple},
%which means that $I$ is integral over a complete intersection
%ideal $J$. For any such minimal reduction $J$ of $I$, one has
%%the minimal primes $\mathfrak{p}$ of $J$ is a minimal prime of $I$
%and
%%\lambda(R_{\mathfrak{p}}/J_{\mathfrak{p}}),\] where $\mathfrak{p}$
%is an arbitrary minimal prime of $I$, or equivalently, of $J$.
%This allows us to state:

\begin{Proposition}\label{E0equi} Let $R$ be a local Cohen--Macaulay
ring with infinite residue field and let $I$ be an equimultiple
ideal of positive height. \begin{itemize} \item[$($a$)$]
$E_0(I)=E_0(J)= \deg(R/J)$ for every minimal reduction $J$ of $I$.
\item[$($b$)$] The ideal $I$ is a complete intersection if and
only if $E_0(I)=\deg(R/I)$ if and only if $E_1(I)=0$.
\end{itemize}
\end{Proposition}
\demo Notice that $J$ is a complete intersection. Furthermore the
minimal primes of $I$ and of $J$ coincide, and hence all have
maximal dimension. For any such prime $\mathfrak{p}$,
$e_0(I_{\mathfrak{p}})=e_0(J_{\mathfrak{p}})=
\lambda(R_{\mathfrak{p}}/J_{\mathfrak{p}})$, proving (a).
% For any minimal prime $\mathfrak p$ of $I$,
Moreover $e_0(I_{\mathfrak{p}})\geq
\lambda(R_{\mathfrak{p}}/I_{\mathfrak{p}})$ and
$e_1(I_{\mathfrak{p}})\geq 0$; either inequality is an equality if
and only if $I_{\mathfrak{p}}$ is a complete intersection
(\cite[Theorem 1]{No}). According to \cite[Theorem]{CN}, the last
condition holds for every $\mathfrak{p}$ if and only if $I$ is a
complete intersection. This proves part (b).\QED

%As to (b), notice that $J$ is a complete intersection and that the
%minimal primes of $I$ and of $J$ coincide. For any such prime
%$\mathfrak{p}$,
%$e_0(I_{\mathfrak{p}})=e_0(J_{\mathfrak{p}})=\lambda(R_{\mathfrak{p}}/J_{\mathfrak{p}})
%$.\QED

%The coefficients $E_i$ share many properties with the ordinary
%Hilbert coefficients, but differ in subtle ways particularly in
%respect to exact sequences. In one case of interest the behavior
%is the expected one. One point that must be kept in sight is that
%for each dimension $s\leq \dim R$ there is a function
%$E_0(\cdot)$, and therefore while considering different modules in
%an exact sequence one must be explicit about which function is
%being used on which module.

\bigskip

The version of Theorem~\ref{e1versusb} for equimultiple ideals can now
be stated. In its proof we will only discuss the points that require a
new justification.

\begin{Theorem}\label{e1versusb2}
Let $R$ be an analytically unramified local Cohen--Macaulay ring
with infinite residue field and let $I$ be an equimultiple ideal
of positive height. Let $A$ and $B$ be distinct graded
$R$--subalgebras of $R[t]$ with
\[ R[It] \subset A \subsetneq B \subset \overline{R[It]}\]
and assume that $A$ satisfies the condition $S_2$ of Serre.
Then
\[ 0 \leq E_1(I) \leq E_1(A) <  E_1(B) \leq \overline{E}_1(I) \leq
 b(I)E_0(I).\]
\end{Theorem}
\demo Let $g$ be the height of $I$ and $\p$ a minimal prime of
 $I$. By Theorem~\ref{e1versusb}, $e_1(A_{\p})
 \leq e_1(B_{\p}) $ and
$e_1(A_{\p})=e_1(B_{\p})$ only when $A_{\p}=B_{\p}$. Now
$A_{\p}=B_{\p}$ for every minimal prime ${\p}$ of $I$ is
equivalent to saying that the $R$--annihilator $L$ of $C=B/A$ is
an ideal of height at least $g+1$. Since $I$ is equimultiple of
height $g$, we conclude that the height of the $A$--ideal $LA$ is
at least $2$. As $LAB \subset A$ and $A$ satisfies the condition
$S_2$, it would follow that $A=B$. This proves the asserted
inequalities except for the last one. To see the last inequality
notice that $b(J_{\mathfrak{p}}) \leq b(I)$ for every minimal
reduction $J$ of $I$, and apply Remark~\ref{REMARK1}.
 \QED

\begin{Remark}{\rm The proof of Theorem 4.2
%~\ref{elversusb2}
shows that when passing from the algebra $A$ to $B$, one of the
values $e_1(A_{\mathfrak{p}})$ increases. Thus, the integer
$\sum_{\mathfrak{p}}e_1(A_{\mathfrak{p}})$ would give tighter
control. Padding the summands with the $\deg(R/\mathfrak{p})$ into
an `Ersatzintegral' however provides a value that becomes
`visible', unlike the $e_1(A_{\mathfrak{p}})$. }\end{Remark}

It is also possible to derive sharper estimates for equimultiple
ideals based on the bounds of Theorem~\ref{maintheorem}. We will
indicate some of these by making use of a very general inequality
for $E_0(I)$ that arises from Lech's formula (\cite[Theorem
3]{Lech}).

\begin{Proposition} Let $R$ be an equidimensional and catenary local Nagata
ring and let $I$ be an ideal of height $g$. Then
\[ E_0(I) \leq g!\,\deg(R/\overline{I})\deg(R).\]
\end{Proposition}

\demo We estimate $E_0(I)$ as given above using Lech's inequality.
Indeed, adding over all minimal primes $\mathfrak p$ of $I$ of
height $g$ we obtain
\begin{eqnarray*}
 E_0(I) &=& \sum e_0(I_{\mathfrak p})\deg(R/{\mathfrak p})\\
  &=& \sum e_0(\overline{I}_{\mathfrak
                p})\deg(R/{\mathfrak p})\\
        &\leq & \sum g! \, \lambda(R_{\mathfrak p}/\overline{I}_{\mathfrak
                p})\deg(R_{\mathfrak p})\deg(R/{\mathfrak p})\\
%&=& g!\,\bigl[\sum  \lambda(R_{\mathfrak
%p}/\overline{I}_{\mathfrak
%                p})\deg(R/{\mathfrak p})\bigr]\deg(R_{\mathfrak p})\\
&\leq & g!\,\bigl[\sum  \lambda(R_{\mathfrak
p}/\overline{I}_{\mathfrak
                p})\deg(R/{\mathfrak p})\bigr]\deg(R)\\
&=& g!\,\deg(R/\overline{I})\deg(R),
\end{eqnarray*}
where we have used the fact that $\deg(R_{\mathfrak p})\leq
\deg(R)$ by our assumption on $R$ (\cite[40.1]{Nag}). \QED

%\lambda(R_{\mathfrak p}/I^n_{\mathfrak p})\deg(R/

%{\mathfrak p}

\bigskip

\begin{Theorem}\label{MT} Let $R$ be a reduced
local Cohen--Macaulay ring and let $I$ be an equimultiple ideal
of height $g >0$.
\begin{itemize}
\item[$($a$)$]
If in addition $R$ is an algebra essentially of finite type over a
perfect field $k$ with type $t$, and $\delta \in {\rm Jac}_k(R)$
is a non zerodivisor, then
\begin{eqnarray*}
\overline{E}_1(I)  & \leq & \frac{t}{t + 1} \ \bigl[ (g-1)E_0(I) +
E_0(I +\delta R/\delta R)\bigr] \quad and \\ \overline{E}_1(I)  &
\leq &
 \frac{t}{t + 1} \ \bigl[ (g-1) \, g! \, \deg(R/\overline{I})\, \deg(R) +
 (g-1)!\, \deg(R/(\overline{I + \delta R}))\,\deg(R/\delta R) \bigr] .
\end{eqnarray*}
\item[$($b$)$]
If the assumptions of $($a$)$ hold, then
\begin{eqnarray*}
\overline{E}_1(I)  & \leq & (g-1)\,\bigl[E_0(I)
-\deg(R/\overline{I})\bigr] + E_0(I +\delta R/\delta R) \quad and \\
\overline{E}_1(I)  & \leq & (g-1)\,\deg(R/\overline{I})\, \bigl[
\, g! \, \deg(R) - 1\bigr]
 %\\ &  &
 +  (g-1)! \, \deg(R/(\overline{I + \delta
R}))\,\deg(R/\delta R).
\end{eqnarray*}
\item[$($c$)$] If $R$ is Nagata and $R/\m$ is infinite, then
\begin{eqnarray*}
\overline{E}_1(I) & \leq & b(I) \ {\rm min} \, \{\frac{t}{t + 1}
\, E_0(I),
E_0(I) -\deg(R/\overline{I})\} \\
&\leq & b(I) \ \frac{t}{t + 1} \ g! \,\deg(R/\overline{I})
\,\deg(R).
\end{eqnarray*}
\end{itemize}
\end{Theorem}
\demo
 To estimate $\overline{E}_1(I)$ in the equimultiple case we start
 from
\[ \overline{E}_1(I)= \sum_{\mathfrak{p}}
\overline{e}_1(I_{\mathfrak{p}})\deg(R/\mathfrak{p}) \] with the
sum taken over the minimal primes $\mathfrak {p}$ of $I$, and make
use of Theorem~\ref{maintheorem} to bound the
$\overline{e}_1(I_{\mathfrak{p}})$ in  terms of the
${e}_0(I_{\mathfrak{p}})$ and ${e}_0((I+ \delta R/ \delta
R)_{\mathfrak{p}})$. Notice that
%if $\p$ is a minimal prime of $I$,
either $\delta \not\in \p$ and then ${e}_0((I+ \delta
R/\delta R)_{\mathfrak{p}})=0$, or $\delta \in \p$ and $\p/\delta
R$ is also a minimal prime of $I+ \delta R/\delta R$, in which
case $I+ \delta R/\delta R$ has height $g-1$. We now process the
summation as in the proof of Proposition 4.4 for the two ideals of
two different rings. For example in (a) we have
\begin{eqnarray*}
\overline{E}_1(I) &=& \sum_{\mathfrak{p}}
\overline{e}_1(I_{\mathfrak{p}})\deg(R/\mathfrak{p}) \\
& \leq & \sum_{\mathfrak{p}} \frac{{\rm type}(R_{\p})}{{\rm
type}(R_{\p})+ 1}\ \bigl[ (g-1)e_0(I_{\p}) + {e}_0((I+ \delta R/
\delta R)_{\mathfrak{p}}) \bigr] \, \deg(R/\mathfrak{p})
\\ & \leq &
 \frac{t}{t + 1}\ \bigl[ (g-1) \sum_{\mathfrak{p}} e_0(I_{\p})\deg(R/\mathfrak{p})+
  \sum_{\mathfrak{p}} {e}_0((I+ \delta R/ \delta
R)_{\mathfrak{p}})\deg(R/\mathfrak{p})\bigr].
% & = & \frac{t}{t + 1}\ [ (g-1)E_0(I) + E_0(I
%+\delta R/\delta R)].
\end{eqnarray*}
In the last expression the first sum equals $E_0(I)$ and the
second sum is either $0$ or else $E_0(I +\delta R/\delta R)$, in
which case $I+ \delta R/\delta R$ has height $g-1$. We now use
Proposition 4.4 to conclude the proof of part (a).
%(Observe also that in the localization the type may
%increase but then the value $t/1+t$ will decrease, so we are safe when using the original
%value. Now someone pick one of the bounds and go through the
%calculation--please pick a pretty one!
\QED


\bigskip
\bigskip






\begin{thebibliography}{99}

\bigskip

\bibitem{AHu}{I. M. Aberbach and C. Huneke, An improved
Brian\c{c}on--Skoda theorem with applications to the
Cohen--Macaulayness of Rees algebras, {\em  Math. Ann.} {\bf 297}
(1993), 343--369.}

\bibitem{AB}{M. Auslander and D. Buchsbaum, On ramification
theory in noetherian rings, {\em Amer. J. Math.} {\bf 81} (1959),
749--765.}

\bibitem{BrH}{W. Brown and J. Herzog,  One--dimensional local
      rings of maximal and almost maximal length, {\em J. Algebra}
{\bf 151} (1992), 332--347.}

\bibitem{CN}{R. C. Cowsik and M. V. Nori, On the fibres of blowing up,
{\em J. Indian Math. Soc.} {\bf 40} (1976), 217--222.}

\bibitem{D}{D. Delfino,  On the inequality $\lambda(\overline
    R/R)\leq t(R)\lambda(R/{\mathcal C})$ for one--dimensional local
    rings,  {\em J. Algebra} {\bf 169} (1994), 332--342.}


\bibitem{E}{J. Elias, On the deep structure of the blowing--up
of curve singularities, {\em Math. Proc. Camb. Philos. Soc.} {\bf
131} (2001), 227--240.}

\bibitem{E2}{J. Elias, On the first normalized Hilbert coefficient, to appear in {\em J. Pure Appl. Alg.}}

\bibitem{GS}{S. Goto and Y. Shimoda, On the Rees algebras of Cohen--Macaulay
local rings, in {\em Commutative Algebra}, Lect. Notes in Pure and
Appl. Math. {\bf 68}, Marcel Dekker, New York, 1982, pp.
201--231.}

\bibitem{HK}{J. Herzog and E. Kunz, {\em Der kanonische Modul eines Cohen--Macaulay Rings},
Lect. Notes in Math. {\bf 238}, Springer, Berlin--New York, 1971.}


\bibitem{H}{M. Hochster, Presentation depth and the Lipman--Sathaye
Jacobian theorem, in {\em The Roos Festschrift}, Vol. 2, Homology
Homotopy Appl. {\bf 4} (2002),  295--314.}

\bibitem{HH}{M. Hochster and C. Huneke, Tight closure in equal characteristic zero, preprint.}


\bibitem{HM}{S. Huckaba and T. Marley, Hilbert coefficients
and the depths of associated graded rings, {\em J. London Math.
Soc.} {\bf 56} (1997), 64--76.}

\bibitem{Hu2}{C. Huneke, Hilbert functions and symbolic powers,
    {\em Michigan
Math. J.} {\bf 34} (1987), 293--318.}

%\bibitem{Hu}{C. Huneke, Uniform bounds in Noetherian rings,
%{\em Invent. Math.} {\bf 107} (1992), 203--223.}

\bibitem{Itoh}{S. Itoh, Coefficients of normal Hilbert polynomials,
{\em J. Algebra} {\bf 150} (1992), 101--117.}

\bibitem{Lech}{C. Lech, Note on multiplicities of ideals,
{\em Ark. Mat.} {\bf 4} (1960), 63--86.     }

\bibitem{LipmanSathaye}{J. Lipman and A.  Sathaye,  Jacobian ideals and
a theorem of Brian\c{c}on--Skoda, {\em Michigan Math. J.} {\bf 28}
(1981), 199--222.}

\bibitem{LiTe}{J. Lipman and B. Teissier, Pseudorational local rings
and a theorem of Brian\c{c}on--Skoda about integral closures of
ideals, {\em Michigan Math. J.} {\bf 28} (1981),  97--116.}


\bibitem{Nag}{M. Nagata,  {\em Local Rings},  Interscience, New York, 1962.}

\bibitem{N}{E. Noether, Idealdifferentiation und Differente,
{\em J. reine angew. Math.} {\bf 188} (1950), 1--21.}

\bibitem{No}{D. G. Northcott, A note on the coefficients of the
abstract Hilbert function, {\em J. London Math. Soc.} {\bf 35}
(1960), 209--214.}

\bibitem{O}{A. Ooishi, $\delta$--genera and sectional genera of
commutative rings, {\em Hiroshima Math. J.} {\bf 17} (1987),
361--372.}

\bibitem{R}{D. Rees, A note on analytically unramified local
    rings, {\em J. London Math. Soc.}  {\bf 36} (1961),
24--28.}

%\bibitem{ram2}{A. Simis, B. Ulrich and W. V. Vasconcelos,
%Codimension, multiplicity and  integral extensions,
%{\em Math. Proc. Camb. Phil. Soc.}  {\bf 130} (2001), 237--257.}

%\bibitem{emb}{B. Ulrich and W. V. Vasconcelos, On the complexity of
%the integral closure,  {\em Trans. Amer. Math. Soc.}, to appear.}

%\bibitem{VallaB}{G. Valla, {Problems and results on Hilbert functions
%of graded algebras}, in
%{\em Six Lectures on Commutative Algebra},
%Progress in Mathematics {\bf 166},  Birkh\"{a}user, Boston, 1998,
%293--344.}


\bibitem{RV}{M. E. Rossi and G. Valla, The Hilbert function of
the Ratliff--Rush filtration, to appear in {\em J. Pure Appl.
Alg.}}

\bibitem{valval}{P. Valabrega and G. Valla, Form rings and
regular sequences, {\em Nagoya Math.
J.} {\bf 72} (1978), 93--101.}

\bibitem{clos}{W. V. Vasconcelos,  Computing the integral closure of
an affine domain,  {\em Proc. Amer. Math. Soc.} {\bf 113} (1991),
633--638.}

%\bibitem{reds}{W. V. Vasconcelos, The reduction numbers of an ideal,
%{\em J. Algebra} {\bf 216} (1999), 652--664.}

%\bibitem{compu}{W. V. Vasconcelos,
%{\em Computational Methods in Commutative
%Algebra and Algebraic Geometry},
%Springer, Heidelberg, 1998.}

\bibitem{teneadd}{W. V. Vasconcelos, Divisorial extensions and the
computation of integral closures, {\em J. Symbolic Comput.} {\bf
30} (2000), 595--604.}

\bibitem{Z}{O. Zariski and P. Samuel, {\em Commutative Algebra},
Vol. 2, Van Nostrand, Princeton, 1960.}


\end{thebibliography}
\end{document}